\documentclass[12pt]{article}
\usepackage{amsmath,eufrak}
\begin{document}
\title{On the remarkable properties of the pentagonal numbers
\footnote{Delivered to the St.--Petersburg Academy September 4, 1775.
Originally published as
{\em De mirabilis proprietatibus numerorum pentagonalium},
Acta Academiae Scientarum Imperialis Petropolitinae
\textbf{4} (1783), no. 1,
56--75, and reprinted in \emph{Leonhard Euler, Opera Omnia}, Series 1:
Opera mathematica, Volume 3, Birkh\"auser, 1992.
A copy of the original text is available
electronically at the Euler Archive, at http://www.eulerarchive.org. This paper
is E542 in the Enestr\"om index.}}
\author{Leonhard Euler\footnote{Date of translation: May 17, 2005.
Translated from the Latin by Jordan
Bell, 3rd year undergraduate in Honours Mathematics, School of Mathematics
and Statistics, Carleton University, Ottawa, Ontario, Canada. Email:
jbell3@connect.carleton.ca. This translation was written during an NSERC
USRA supervised by Dr. B. Stevens.}}
\date{}

\maketitle

1. To the class of pentagonal numbers I refer not only to those which
are generally and indeed appropriately called so by custom and comprised
of the form $\frac{3nn-n}{2}$, but also those which are made in this
form: $\frac{3nn+n}{2}$. Thus, the general form of all these
numbers will be $\frac{3nn \pm n}{2}$, from which if in place of $n$
is written the succession of numbers 0, 1, 2, 3, 4, etc., 
the following sequence of paired numbers arises:

\begin{tabular}{p{3cm}|p{0.75cm}p{0.75cm}p{0.75cm}p{0.75cm}p{0.75cm}p{0.75cm}p{0.75cm}}
$n$&0,&1,&2,&3,&4,&5,&6\\
\hline
pentagonal&0,&1,&5,&12,&22,&35,&51\\
numbers&0,&2,&7,&15,&26,&40,&57
\end{tabular}

That is, any number taken by $n$ produces two numbers, 
which are written as one above the other, such that the top sequence
contains the pentagonal numbers properly called,
and the bottom indeed those which I refer to by the very
same class, and which come forth after the above sequence.

2. But if we put these numbers in a single sequence in order of their magnitudes,
this progression will unfold:
\[
0,1,2,5,7,12,15,22,26,35,40,51,57,\textrm{ etc.,}
\]
whose order is clearly broken up, seeing that the progression of differences
is
\[
1,1,3,2,5,3,7,4,9,5,11,\textrm{ etc.},\]
which is made from the sequence of natural numbers with others mixed
in. Moreover, it is possible for this sequence to be
led into a continuation, if after every third term a certain
fraction is interpolated. Namely, between the terms 2 and 5 is placed
$\frac{10}{3}$, then indeed $\frac{28}{3}$ between 7 and 12,
similarly $\frac{55}{3}$ between 15 and 22, such that the sequence
will be consummated as:
\[
1,2,\frac{10}{3},5,7,\frac{28}{3},12,15,\frac{55}{3},22,26,\textrm{ etc.};
\]
it follows by continuing along with this rationale that the sequence of differences will
then be:
\[
1,\frac{4}{3},\frac{5}{3},2,\frac{7}{3},\frac{8}{3},3,\frac{10}{3},
\frac{11}{3},4\textrm{ etc.}
\]
Moreover, it is clear that the same sequence is seen if all the triangular
numbers are divided by 3. From this therefore, a beautiful property of
our pentagonal numbers shows itself, seeing that every third shows
itself to be a triangular number.

3. So extraordinary are the properties which can be derived from the general
formula and which moreover are able to be applied to other polygonal
numbers, toward which consequently I do not gaze. At this moment,
it is preferable to recall a proposition of mine
on certain astonishing properties which the
pentagonal numbers above all the other polygonal numbers are gifted with.
Here, a significant property of these numbers comes to mind which I have
already shown, for there to be a very close connection between this
sequence of pentagonal numbers and the progression formed by the sums of
the divisors of the natural numbers, so that by means of this, in fact,
a rule for this most irregular sequence can be assigned, with respect
to which it will be worthwhile to quickly repeat the work. 

4. If for any number $N$ its divisors are collected into a single sum,
then we will signify this sum with the character $\int N$:
then from the natural numbers is born a sequence which when first looked at
is most irregular:

\begin{tabular}{r|ccccccccccc}
$N$&1,&2,&3,&4,&5,&6,&7,&8,&9,&10,&11\\
\hline
$\int N$&1,&3,&4,&7,&6,&12,&8,&15,&13,&18,&12
\end{tabular}

in which the terms 
progress so inordinately, at one time increasing and at another time decreasing,
that any rule for them will be uncovered with difficulty, seeing that
 the sequence of prime numbers themselves is obviously enveloped
in them.

5. I have nevertheless however demonstrated that this progression,
although certainly irregular, can be considered in the class
of recurrent sequences, where each term immediately next
is able to be determined by a certain rule from the preceeding ones.
For, if $\int N$ denotes the sum of all the divisors of the number $N$, with
itself not excepted, it will always be found to be:
\begin{align*}
\int N= \int (N-1) + \int (N-2) - \int (N-5) - \int (N-7)+
\\
\int (N-12) + \int (N-15) -
 \int (N-22) - \int (N-26) + \textrm{ etc.}
\end{align*}
in which the numbers which are successively subtracted from $N$ clearly
make up our sequence of pentagonal numbers
\[
1,2,5,7,12,15,22,26,35,40,\textrm{ etc.}
\]
where the terms that come forth for $n$ taken as odd
have the sign $+$, while indeed when it is even the sign
$-$ is formed.
Then indeed, after each case it is required
for the expressions to be continued thusly, as long as the number written
after
the $\int$ sign  
does not come out to be negative; if the expression $\int (N-N)$ occurs, in place
of it should be written the number $N$. Thus if we select $N=12$, it will
be:
\[
\int 12 = \int 11 + \int 10 - \int 7 - \int 5 + \int 0
\]
and it will therefore be:
\[
\int 12 = 12 + 18 - 8 - 6 + 12=28.
\]
To be sure, if we take up $N=13$, it will be:
\[
\int 13 = \int 12 + \int 11 - \int 8 - \int 6 + \int 1,
\]
that is, it will be:
\[
\int 13 = 28 + 12 - 15 - 12 + 1 = 14.
\]

6. Seeing then that the order in which the sums of divisors procede,
which appears most deservedly as irregular and which
no one has been able to come to a
definite conclusion about, has been able to be explored with the pentagonal
numbers, from which this observation is by all means to be admired highly.
I will set forth moreover another property of this which is indeed closely
connected to what has been related, but which leads to several no less
admirable properties, which will arise equally from the nature of our
pentagonal numbers which has been established.

7. The beginning of all of these remarkable properties is encountered
in the expansion of this infinite product:
\[
S=(1-x)(1-xx)(1-x^3)(1-x^4)(1-x^5)(1-x^6)(1-x^7) (\textrm{ etc.},
\]
for I have demonstrated that if each of the factors are multiplied by the
others in order, then in the end this series results:
\[
S=1-x^1-x^2+x^5+x^7-x^{12}-x^{15}+x^{22}+x^{26}-\textrm{ etc.},
\]
where the exponents for each $x$ form our sequence of pentagonal numbers,
with the rule for the signs $+$ and $-$ that they alternate in pairs
of two, so that the exponents which procede with $n$ taken as
even
have the sign $+$, and indeed the the others which come from it odd
the sign $-$. 
This then deserves our admiration no less than the properties mentioned
before, with no fixed rule apparant from which any connection can be
understood between the expansion of this product and our pentagonal
numbers.

8. Therefore with this series of powers of $x$ equal to the above
infinite product, if we take it equal to nothing, so that
we would have the equation:
\[
0=1-x^1-x^2+x^5+x^7-x^{12}-x^{15}+x^{22}+x^{26}-\textrm{ etc.}
\]
it will involve each of the roots which the former
product forces to be equal to nothing. According to the first factor
$1-x$ it will be of course $x=1$; from the second factor $1-xx$ it will
be both $x=+1$ and $x=-1$; from the third factor $1-x^3$ originate
these three roots:
\[
1^{\circ}) \, x=1, \qquad 2^{\circ}) \, x= -\frac{1+\surd{-3}}{2}, \qquad
3^{\circ}) \, x=-\frac{1-\surd{-3}}{2};
\]
then again, from the fourth factor $1-x^4=0$ emerge these four roots:
\[
1^{\circ}) \, x=+1, \qquad 2^{\circ}) \, x=-1, \qquad
3^{\circ}) \, x=+\surd{-1} \qquad \textrm{ and } \qquad 4^{\circ}) \, x=-\surd{-1}; 
\]
next the fifth factor $1-x^5=0$ yields these five roots:
\begin{align*}
1^{\circ}) \, x=1, \qquad 2^{\circ})
\, x=\frac{-1-\surd{5}+\surd(-10+2\surd{5})}{4}, \qquad \\
3^{\circ}) \, x= \frac{-1-\surd{5}-\surd(-10+2\surd{5})}{4}, \qquad
4^{\circ}) \, x= \frac{-1+\surd{5}+\surd(-10-2\surd{5})}{4}, \qquad
\\
5^{\circ}) \, x= \frac{-1+\surd{5}-\surd(-10-2\surd{5})}{4};
\end{align*}
in turn the sixth factor produces these six roots:
\begin{align*}
1^{\circ}) \, x=1, \qquad 2^{\circ}) \, x=-1, \qquad
3^{\circ}) \, x=\frac{+1+\surd{-3}}{2},\\
4^{\circ}) \, x=\frac{+1-\surd{-3}}{2}, \qquad
5^{\circ}) \, x=\frac{-1+\surd{-3}}{2},\\
6^{\circ}) \, x=\frac{-1-\surd{-3}}{2} \qquad \textrm{ etc. etc.}
\end{align*}

9. Consequently then, it is clear for all the roots of unity to
simultaneously be roots of our equation, and if we consider the general
problem 
in which it is to be put
$1-x^n=0$, it is at once apparent for $x=1$
to always be a root, and if $n$ is an even number for another root to
be $x=-1$. For the other roots, the trionomial factors of $1-x^n$
ought to be considered, which, making use of what has already been set
forth elsewhere, are contained in the general form:
\[
1-2x\cos{\frac{2i\pi}{n}}+xx,
\]
taking successively for $i$ all the integral numbers no greater than
$\frac{1}{2}n$. 
By equating this factor with nothing, these two roots are elicited:
\begin{align*}
x=\cos{\frac{2i\pi}{n}}+\surd{-1}\sin{\frac{2i\pi}{n}} \textrm{ and}\\
x=\cos{\frac{2i\pi}{n}}-\surd{-1}\sin{\frac{2i\pi}{n}}.
\end{align*}
For then it will be in turn:
\[
x^n=\cos{2i\pi} \pm \surd{-1} \sin{2i\pi}.
\]
Moreover, when it is $\cos{2i\pi}=1$ and $\sin{2i\pi}=0$, then
$x^n=1$, from which, if for $n$ and $i$ are taken all the successive
integral numbers, this formula:
\[
x=\cos{\frac{2i\pi}{n}} \pm \surd{-1}\sin{\frac{2i\pi}{n}}
\]
will produce all the roots of our equation
\[
0=1-x-x^2+x^5+x^7-x^{12}-x^{15}+x^{22}+x^{26}-\textrm{etc.}
\]
so that plainly we may have the power to assign all the roots of this equation.

10. But if we indicate all the roots of this equation by the letters
$\alpha, \beta,\gamma\, \delta, \epsilon$, etc., the factors
will be $\frac{1-x}{\alpha}, \frac{1-x}{\beta}, \frac{1-x}{\gamma},
\frac{1-x}{\delta}$, etc., from which we deduce from the nature of
equations for the sum of all these fractions $\frac{1}{\alpha}+
\frac{1}{\beta}+\frac{1}{\gamma}+\frac{1}{\delta}+\textrm{ etc.}=1$,
then in fact for the sum of the products two apiece to be equal to -1, 
the sum of products three apiece to indeed be equal to 0, the sum of the products
four apiece to be equal to 0, the sum of the products five apiece to be equal
to -1, the sum of the products six apiece to be equal to -1, etc. Then
on the other hand, we conclude as before for the sum of the squares
of these fractions to be
\[
\frac{1}{\alpha^2}+\frac{1}{\beta^2}+\frac{1}{\gamma^2}+\frac{1}{\delta^2}+
\textrm{ etc.}=3,
\]
the sum of the cubes
\[
\frac{1}{\alpha^3}+\frac{1}{\beta^3}+\frac{1}{\gamma^3}+\frac{1}{\delta^3}+
\textrm{ etc.}=4,
\]
the sum of the biquadratics
\[
\frac{1}{\alpha^4}+\frac{1}{\beta^4}+\frac{1}{\gamma^4}+\frac{1}{\delta^4}+
\textrm{ etc.}=7,
\]
and thusly as before; in this certainly no order is observed.

11. Insofar as we have the fractions $\frac{1}{\alpha},\frac{1}{\beta},
\frac{1}{\gamma}$, etc., moreover furthermore $\alpha,\beta,\gamma$, etc.
prevail as roots themselves. For if $\alpha$ is a root
of our equation, we can express this root in the form:
\[
\cos{\frac{2i\pi}{n}} \pm \surd{-1}\sin{\frac{2i\pi}{n}}.
\]
Then moreover it is
\[
\frac{1}{\alpha}=\frac{1}{\cos{\frac{2i\pi}{n}} \pm \surd{-1}
\sin{\frac{ 2i\pi}{n}}}=\cos{\frac{2i\pi}{n}} \mp \surd{-1}
\sin{\frac{2i\pi}{n}},
\]
which will likewise be a root of our equation; from this it is apparent
that if $\frac{1}{\alpha}$ is a root of our equation, then $\alpha$
will also be a root.

12. Therefore, were $\alpha$ to denote any root of the equation
$1-x^n=0$, seeing that then it will also be a root of our
equation
\[
1-x-xx+x^5+x^7-x^{12}-x^{15}+\textrm{ etc.}=0,
\]
then therefore it will be $\alpha^n=1$. In addition, in fact all
powers of $\alpha$ will also be roots of the equation $1-x^n=0$.
If for in place of $x$ we write $\alpha \alpha$, it will be $1-x^n=1-\alpha^{2n}$.
With it moreover $\alpha^n=1$, it is clear then for it to be
$\alpha^{2n}=1$, and therefore $1-\alpha^{2n}=0$, with respect to
which the same is clear for the cube $\alpha^3$ and all higher powers.
Then therefore it follows that it is:
\[
\alpha^{n+1}=\alpha \textrm{ and } \alpha^{n+2}=\alpha \alpha \textrm{ and }
\alpha^{n+3}=\alpha^3.
\]
It follows that in general it will be $\alpha^{in+\lambda}=\alpha^\lambda$.

13. If then $\alpha$ denotes some arbitrary
root of our equation, so that it is
thus $\alpha^n=1$, if in place of each $x$ we write $\alpha$, clearly
this series arises:
\[
1-\alpha^1-\alpha^2+\alpha^5+\alpha^7-\alpha^{12}-\alpha^{15}+\alpha^{22}+
\textrm{ etc.}=0.
\]
Thereafter indeed, by putting $x=\alpha \alpha$ it will be:
\[
1-\alpha^2-\alpha^4+\alpha^{10}+\alpha^{14}-\alpha^{24}-\alpha^{30}+
\alpha^{44}+\textrm{ etc.}=0,
\]
and in general if in place of $x$ we write $\alpha^i$, denoting with $i$
any integral number, it will further be:
\[
1-\alpha^i-\alpha^{2i}+\alpha^{5i}+\alpha^{7i}-\alpha^{12i}-
\alpha^{15i}+\alpha^{22i}+\textrm{ etc.}=0.
\]
Besides this it occurs on the other hand that when negative numbers are taken
for $i$, as we have already shown, each of the roots will be
$\frac{1}{\alpha^2},\frac{1}{\alpha^3},\frac{1}{\alpha^4},\frac{1}{\alpha^5}$,
etc.

14. Seeing here that
we have assumed $\alpha$ to be a root of the equation
$1-x^n=0$, we have run through the succession of cases in which $n$
is either 1, or 2, or 3, or 4, etc; and indeed in the first, where $n=1$, 
by necessity it is $\alpha=1$. By substituting this value into
our general equation, this form is induced:
\[
1-1-1+1+1-1-1+1+\textrm{ etc.;}
\]
this series is evidently composed from infinitely
many periods of which each
contains these terms: $1-1-1+1$, according to which the value of each
of each of the periods is equal to 0, for which reason then the infinitely
many periods likewise have their sum taken as equal to 0. 
Moreover, it should be understood that by continuing, if now the infinitely
many periods above are run through after one [sic] term
is constituted, the sum will
be equal to 0; if three are constituted, the sum will be will -1, and if
four are constituted it will be equal to 0, in which case an entire
period has been added.  
Therefore, since this infinite sum of numbers never terminates,
the sum of this infinite series will be kept between the four sums
1, 0, -1, 0 just considered. The middle is obtained if these four
numbers are aggregated into one that is then divided by four; then moreover
it is apparent to come out as 0, which therefore is rightly decreed the sum of
our series.

15. One may know that a reckoning similar to this can be carried out,
which is usually referred to by the sum of the Leibniz series
$1 - 1 + 1 - 1 + 1 - 1 + 1 - 1$ + etc. that is equal to $\frac{1}{2}$; with this being
assumed, the truth of the present assertion can shine forth. Seeing that
it is:
\begin{align*}
1-1+1-1+1-1+\textrm{ etc.}=\frac{1}{2}, \textrm{ it will be }\\
-1+1-1+1-1+1-\textrm{ etc.}=-\frac{1}{2},
\end{align*}
hence by combining these two series it will be:
\[
1-1-1+1+1-1-1+1+1-1-\textrm{ etc.}=0.
\]

16. We now consider the case in which $n=2$ and $\alpha \alpha=1$, in
which $\alpha$ is either $+1$ or $-1$. However, we retain the letter
$\alpha$ for designating either of these, and because it is
\[
\alpha^3=\alpha, \alpha^4=1, \alpha^5=\alpha, \alpha^6=1,
\textrm{ etc.},
\]
substiting what has been made into our general equation induces this
form:
\[
1-\alpha-1+\alpha+\alpha-1-\alpha+1|+1-\alpha-1+\alpha+\alpha-1-\alpha+1+
\textrm{ etc.;}
\]
this series is advanced by a certain period, which is continuosly repeated,
and every one of which consists of these eight terms:
\[
1-\alpha-1+\alpha+\alpha-1-\alpha+1,
\]
of which the sum is 0, and however large, the complete period is certain
to disappear. If indeed the first one, two, three or all
the way
the eight terms are began with, the sums will be obtained in the following
way:

\begin{tabular}{p{4cm}|p{3cm}}
given the first terms&the sum will be\\
one&1\\
two&$1-\alpha$\\
three&$-\alpha$\\
four&0\\
five&$\alpha$\\
six&$\alpha-1$\\
seven&-1\\
eight&0 
\end{tabular}

of which eight the sum of their aggregate is 0, from which
we can safely conclude for the sum of the whole of this series that we have made,
continued infinitely, to be equal to 0.

17. From this it is clear for the sum of this periodic series to likewise
be equal to nothing, for any value the letter $\alpha$ may have;
 truly then, consideration is led to of the periods which arise under
the value of $\alpha$ for which $\alpha \alpha=1$, wherefore this series is
able to be separated into two parts, of which one contains only
unities, and the other indeed only the letters $\alpha$. It is necessary
that the sum of both be equal to nothing, so that it is:
\begin{align*}
1-1-1+1, +1-1-1+1, +1-1-1+1, \textrm{ etc.}=0\\
-\alpha+\alpha+\alpha-\alpha,-\alpha+\alpha+\alpha-\alpha,-\alpha+
\alpha+\alpha-\alpha, \textrm{ etc.}=0,
\end{align*}
for which the truth of both is apparent from what has been
established earlier.

18. In a similar way, for the event of having the cube root of 1 by placing
$\alpha^3=1$, and 
seeing that several terms of the period will extend, we will refer to the
general series by two terms written one below the other, so that in general:
\begin{align*}
\left. \begin{array}{r}
1-\alpha+\alpha^5-\alpha^{12}+\alpha^{22}-\alpha^{35} \textrm{ etc.}\\
-\alpha^2+\alpha^7-\alpha^{15}+\alpha^{26}-\alpha^{40} \textrm{ etc.}
\end{array}
\right\}=0.
\end{align*}
Seeing that if it is now supposed $\alpha^3=1$, such that it is
\[
\alpha^4=\alpha, \alpha^5=\alpha^2, \alpha^6=1, \alpha^7=\alpha,
\textrm{ etc.}
\]
the follows periodic progression will extend:
\begin{align*}
\left. \begin{array}{r}
1-\alpha+\alpha^2-1+\alpha-\alpha^2+1|-\alpha+\alpha^2-1+\alpha-\alpha^2+1|\\
-\alpha^2+\alpha-1+\alpha^2-\alpha|+1-\alpha^2+\alpha-1+\alpha^2-\alpha|+1
\end{array} \right. \textrm{ etc.}
\end{align*}
equal to nothing, 
in which each period is constituted by twelve terms of the three general
ones, namely $1, \alpha, \alpha^2$. It is noticed without difficulty
for the general terms taken for this series to yield nothing, for the
unities constitute this series:
\[
1-1-1+1, +1-1-1+1, +1-1-1+1, \textrm{ etc.}=0,
\]
and indeed the letters $\alpha$ and $\alpha \alpha$ constitute the following
series:
\begin{align*}
-\alpha+\alpha+\alpha-\alpha, -\alpha+\alpha+\alpha-\alpha, -\alpha+
\alpha+\alpha-\alpha, \textrm{ etc.}=0\\
-\alpha^2+\alpha^2+\alpha^2-\alpha^2,-\alpha^2+\alpha^2+\alpha^2-\alpha^2,
-\alpha^2+\alpha^2+\alpha^2-\alpha^2,\textrm{ etc.}=0.
\end{align*}
It is clear for the sums of both of these to be equal to nothing.

19. We will further now consider the biquadratic roots of unity, for which it
is $\alpha^4-1$, and the following periodic series is produced:
\begin{align*}
1-\alpha+\alpha-1+\alpha^2-\alpha^3+\alpha^3-\alpha^2+1|-\alpha+\alpha-1+
\textrm{ etc.}\\
-\alpha^2+\alpha^3-\alpha^3+\alpha^2-1+\alpha-\alpha|+1-\alpha^2+\alpha^3-
\alpha^3+\textrm{ etc.}
\end{align*}
for which each period is constituted from 16 terms. The general four
which have been related are shown in the following four
series each to be equal to nothing:
\begin{align*}
1-1-1+1, +1-1-1+1, +1-1-1+1, \textrm{ etc.}=0\\
-\alpha+\alpha+\alpha-\alpha, -\alpha+\alpha+\alpha-\alpha,-\alpha+\alpha+
\alpha-\alpha, \textrm{ etc.}=0\\
-\alpha^2+\alpha^2+\alpha^2-\alpha^2,-\alpha^2+\alpha^2+\alpha^2-\alpha^2,
-\alpha^2+\alpha^2+\alpha^2-\alpha^2, \textrm{ etc.}=0\\
+\alpha^3-\alpha^3-\alpha^3+\alpha^3,+\alpha^3-\alpha^3-\alpha^3+\alpha^3,
+\alpha^3-\alpha^3-\alpha^3+\alpha^3, \textrm{ etc.}=0.
\end{align*}

20. Although now our conclusion about the higher roots has been
well confirmed, it is however necessary to investigate the additional
case of $\alpha^5=1$, seeing that here not all the powers less than five
occur. If then it is $\alpha^5=1$, this periodic series will advance:
\begin{align*}
\left. \begin{array}{r}
1-\alpha+1-\alpha^2+\alpha^2-1+\alpha-1+\alpha^2-\alpha^2+1|-\alpha+1\\
-\alpha^2+\alpha^2-1+\alpha-1+\alpha^2-\alpha^2+1-\alpha|+1-\alpha^2+\alpha^2
\end{array} \right. \textrm{ etc.}
\end{align*}
in which the powers $\alpha^3$ and $\alpha^4$ are thoroughly excluded.
Since each of the periods is constituted by exactly 20 terms, it
is necessary for the other powers to occur frequently; the separate
three that are taken occur in the following periodic series:
\begin{align*}
1+1-1-1-1-1+1+1,1+1-1-1-1-1+1+1,\textrm{ etc.}=0\\
-\alpha+\alpha+\alpha-\alpha, -\alpha+\alpha+\alpha-\alpha,
-\alpha+\alpha+\alpha-\alpha, \textrm{ etc.}=0\\
-\alpha^2+\alpha^2-\alpha^2+\alpha^2-\alpha^2+\alpha^2-\alpha^2+\alpha^2,
-\alpha^2+\alpha^2-\alpha^2+\alpha^2, \textrm{ etc.}=0.
\end{align*}
By now the truth of the series of $\alpha$ themselves is clear
from the preceding; for the remaining two, of which the period consists of
eight terms, if the following which has been established earlier is
examined, even now they will be discovered to be equal to
nothing. For, not only do the terms alone in the periods eliminate
each other,  but also the terms of the basis from which the
series is formed. Then from the series of unities is derived this
base series:
\[
1,2,1,0,-1,-2,-1,0,
\]
of which the sum again vanishes; by the same experience the series
of squares is disposed of.

21. From this it is now abundantly clear that these same properties
will keep their place for higher roots, however many terms each period
is composed of. This is certainly most amazing, with these properties
able to occur in no other series of powers, and is a deep
property of the series of pentagonal numbers.

22. So that moreover we can put eyes on the problem in general, it
will be $\alpha^n=1$, from which are born periods consisting
of $4n$ terms, which will be either 1, or $\alpha$, or $\alpha^2$,
or $\alpha^3$, etc. While for the most part not all the power less than
$\alpha^n$ will occur, the periods of each power of $\alpha$
are typically constituted from more than four terms.  
However, not only do the terms themselves of the period always
destroy each other but also terms of the basis. Therefore if we
consider the powers $\alpha^r$, with $r$ arising as a number less
than $n$, according to all the terms taken from our sequence of
pentagonal numbers which divided by $n$ leave behind the residue $r$.
If these terms are given with the sign put in front, such a series
comes forth:
\[
\pm \alpha^r \pm \alpha^r \pm \alpha^r \pm \alpha^r
\pm \alpha^r \pm \alpha^r \pm \alpha^r \pm \alpha^r
\pm \textrm{ etc.,}
\]
which always consists of a certain period with a rule for the
signs $+$ and $-$, so that of each period all the terms
that are taken at once destroy each other, and moreover the same
for the series that occurs of the basis.

23. Truly, along with these properties that have been related
so far there are countless others no less extraordinary that 
are drawn forth after these. For if $\alpha$ is a root
of an power $n$ of unity, so that $1-\frac{x}{\alpha}$ is a 
factor of the form $1-x^n$, it is evident for it to be a factor
of the forms $1-x^{2n}, 1-x^{3n}, 1-x^{4n}$, etc. into infinity.
From this, because these forms are all factors of our progression:
\[
1-x-xx+x^5+x^7-x^{12}-\textrm{ etc.},
\]
each root $\alpha$ occurs not only once but indeed infinitely in this equation,
 so that this equation has infinitely many roots equal to $\alpha$.

24. We know on the other hand from the nature of equations, that if
some equation
\[
1+Ax+Bxx+Cx^3+Dx^4+\textrm{ etc.}=0
\]
has two root equal to $\alpha$ then in fact $\alpha$ is a root
of the equation made by differentiating, namely:
\[
A+2Bx+3Cxx+4Dx^3+\textrm{ etc.}=0,
\]
and if it has three roots equal to $\alpha$ then 
in addition $\alpha$ will likewise be a root of
the equation that is made by differentiating, after we have
multiplied this
differentiated equation by $x$:
\[
1^2\cdot A + 2^2\cdot Bx + 3^2\cdot Cxx + 4^2\cdot Dx^3 + \textrm{ etc.}=0,
\]
from which if these equation has $\lambda$ equal roots, each
of which are equal to $\alpha$, then it will always be:
\[
1^\lambda \cdot A + 2^\lambda \cdot B\alpha + 3^\lambda \cdot C\alpha \alpha + 
4^\lambda \cdot D\alpha^3 + \textrm{ etc.}=0,
\]
from which  if we multiply uniformly this equation by
$\alpha$, it will then be:
\[
1^\lambda \cdot A\alpha + 2^\lambda \cdot B\alpha^2 + 3^\lambda \cdot C\alpha^3 + 
4^\lambda \cdot D\alpha^4 + \textrm{ etc.}=0.
\]

25. Therefore with it having been put $\alpha^n=1$, our equation formed
from the pentagonal numbers:
\[
 1-x^1-x^2+x^5+x^7-x^{12}-x^{15}+\textrm{ etc.}=0
\]
has infinitely many roots equal to $\alpha$, and it will actually
be that $\alpha$ is a root of all equations contained in this
general form:
\[
-1^\lambda x -2^\lambda x^2 + 5^\lambda x^5 + 7^\lambda x^7 -
12^\lambda x^{12} - \textrm{ etc.}=0,
\]
where $\lambda$ is taken as any integral number. Therefore it will
always be:
\[
-1^\lambda \alpha - 2^\lambda \alpha^2 + 5^\lambda \alpha^5
+ 7^\lambda \alpha^7 - 12^\lambda \alpha^{12}-\textrm{ etc.}=0.
\]

26. In order to clearly illustrate this we take $\alpha=1$, and it will
always be:
\[
-1^\lambda -2^\lambda + 5^\lambda +7^\lambda -12^\lambda -15^\lambda +
\textrm{ etc.}=0,
\]
and for the case of $\lambda=0$, we have already examined the truth
of this equation. Therefore if it were $\lambda=1$, it is to be
shown
for the sum of this sucession which diverges to infinity to be equal to 0:
\[
 -1-2+5+7-12-15+22+26-\textrm{ etc.}
\]
Because this series is interrupted, that is, it has been mixed from two
series, we may consider it separated into two parts by putting:
\begin{align*}
s=-1+5-12+22-35+\textrm{ etc., and}\\
t=-2+7-15+26-40+\textrm{ etc.}
\end{align*}
and it ought to be shown to be $s+t=0$.

27. From the doctrine of series which proceed with their signs alternating,
 just as $A-B+C-D+\textrm{ etc.}$ continues, 
the sum of this series progressing to infinity will be equal to:
\[
\frac{1}{2}A-\frac{1}{4}(B-A)+\frac{1}{8}(C-2B+A)-\frac{1}{16}(D-3C+3B-A)
\textrm{ etc.,}
\]
This useful rule is set forth by differences, namely by a rule
based on the signs. From the series of numbers $A,B,C,D,E$, etc. is formed
the series of differences, so that each term of this series is subtracted
from the following one: it will be $a,b,c,d$, etc. Again, in the very
same, from this series of differences is formed a series of second differences,
which will be $a', b', c', d'$, etc., from which likewise is the series
of third differences, which will be $a'', b'', c'', d'', e''$, etc., and
in this way further differences are reached for until they are unchanging.
Then moreover from the first terms of all of these series the sum of
the proposed series may thus be determined, such that it will be:
\[
\frac{1}{2}A-\frac{1}{4}a+\frac{1}{8}a'-\frac{1}{16}a''+\frac{1}{32}a'''-
\frac{1}{64}a''''+\textrm{ etc.}
\]

28. By this rule it is established, with the signs changed it will be:
\begin{align*}
-s=1-5+12-22+35-51+70-\textrm{ etc., and}\\
-t=2-7+15-26+40-57+77-\textrm{ etc.},
\end{align*}
 whose terms are arranged in the following way with their differences
written underneath them:

\begin{tabular}{c||c}
1,5,12,22,35,51,70 etc.&2,7,15,26,40,57,77 etc.\\
4,7,10,13,16,19&5,8,11,14,17,20\\
3,3,3,3,3&3,3,3,3,3\\
0,0,0,0&0,0,0,0
\end{tabular}

Therefore from this it is gathered to be:
\begin{align*}
-s=\frac{1}{2}-\frac{4}{4}+\frac{3}{8}=-\frac{1}{8}, \textrm{ that is, }
s=\frac{1}{8}, \textrm{ and moreover,}\\
-t=\frac{2}{2}-\frac{5}{4}+\frac{3}{8}=\frac{1}{8}, \textrm{ that is, }
t=-\frac{1}{8},
\end{align*}
from which it is clearly combined to be $s+t=0$.

29. Although the rules which these properties are supported by clearly
leave no doubt, it will be by no means useless to exhibit
further the truth for the case $\lambda=2$, for it to be:
\[
-1^2-2^2+5^2+7^2-12^2-15^2+22^2+\textrm{ etc.}=0.
\]
In the same way as before this series is divided into two, which with
their signs changed are:
\begin{align*}
s=+1^2-5^2+12^2-22^2+35^2-51^2+\textrm{ etc.}\\
t=2^2-7^2+15^2-26^2+40^2-57^2+\textrm{ etc.},
\end{align*}
and in order to find the sum of the first, the following procedure is
carried out:

\begin{tabular}{p{4cm}p{9cm}}
\textrm{Series}&1, 25, 144, 484, 1225, 2601, 4900\\
\textrm{Diff. I}&24, 119, 340, 741, 1376, 2299\\
\textrm{Diff. II}&95, 221, 401, 635, 923\\
\textrm{Diff. III}&126, 180, 234, 288\\
\textrm{Diff. IV}&54, 54, 54\\
\textrm{Diff. V}&0, 0
\end{tabular}

Hence it will then be:
\[
s=\frac{1}{2}-\frac{24}{4}+\frac{95}{8}-\frac{126}{16}+\frac{54}{32}=
+\frac{3}{16}.
\]
In the same way for the other series:

\begin{tabular}{p{4cm}p{9cm}}
\textrm{Series}&4, 49, 225, 676, 1600, 3249, 5929\\
\textrm{Diff. I}&45, 176, 451, 924, 1649, 2680\\
\textrm{Diff. II}&131, 275, 473, 725, 1031\\
\textrm{Diff. III}&144, 198, 252, 306\\
\textrm{Diff. IV}&54, 54, 54\\
\textrm{Diff. V}&0, 0
\end{tabular}

Then it is concluded:
\[
t=\frac{4}{2}-\frac{45}{4}+\frac{131}{8}-\frac{144}{16}+\frac{54}{32}=
-\frac{3}{16}.
\]
On account of this, it prevails for the total sum to be $s+t=0$.

30. Furthermore, we now consider the square root, i.e. when it is
$\alpha^2=1$, and hence this series arises:
\[
-1^\lambda \cdot \alpha - 2^\lambda + 5^\lambda \cdot \alpha+
7^\lambda \alpha - 12^\lambda  - 15^\lambda \cdot \alpha + 22^\lambda+
26^\lambda - \textrm{ etc.}=0,
\]
from which if we separate the terms comprised of unity and $\alpha$
from each other, we will obtain two series both equal to nothing, namely:
\[
-2^\lambda - 12^\lambda + 22^\lambda + 26^\lambda - 40^\lambda - 70^\lambda
+ 92^\lambda + \textrm{ etc.}=0
\]
and
\[
-1^\lambda \cdot \alpha + 5^\lambda \cdot \alpha + 7^\lambda \cdot \alpha -
15^\lambda \cdot \alpha - 35^\lambda \cdot \alpha + 51^\lambda \cdot \alpha
+ 57^\lambda \cdot \alpha - \textrm{ etc.}=0.
\]
Insofar as we may want to display the truth of these series in the same
manner which we have made use of before, each one ought to be broken
into four others, until in the end we reach unchanging differences. Truly,
if this work were to have been taken up, it will be certain for the
aggregate of the sums of all the parts to be equal to 0.

31. Now most generally the whole matter is embraced, where it will
be $\alpha^n=1$, and we search for a series which
contains all powers $\alpha^r$. We select to this end from all
our pentagonal numbers which when divided by $n$ leave a residue
the same as $r$. Then if the pentagonal numbers are $A,B,C,D,E$, etc., all
 of course of the form $\gamma n+r$ and of which the sign $\pm$,
which is appropriate for them, is observed carefully. Then
indeed it will always be:
\[
\pm A^\lambda \pm B^\lambda \pm C^\lambda \pm D^\lambda \pm \textrm{ etc.}=0,
\]
for which the exponent $\lambda$ may take any integral value. Then in this
most general form, all series which we have elicited so far, and of which we have
shown for the sums to be equal to nothing, are contained.

\end{document}